\documentclass{birkjour}
\usepackage{amsxtra}
\usepackage{amsthm,mathrsfs}
\usepackage{amsmath,amsthm,amssymb,amsrefs,bbm,color,enumitem,esint,esvect,float,graphicx,mathrsfs, hyperref}
\begin{document}
\begin{abstract}
In this paper we offer alternate upper bound for the operator $\Pi_b^*\Pi_d$ to the ones present in literature, thus extending the known upper bounds from the $L^2(\mathbb{R})$ setting to $L^p(w)$, for $1<p<\infty,$ and a Muckenhoupt weight $w$. In the $L^2(w)$ setting, we fully characterize the boundedness of the operator. 
\end{abstract}
 \newtheorem{theorem}{Theorem}[section]
 \newtheorem{corollary}[theorem]{Corollary}
 \newtheorem{lemma}[theorem]{Lemma}
 \newtheorem{proposition}[theorem]{Proposition}
 \theoremstyle{definition}
 \newtheorem{definition}[theorem]{Definition}
 \theoremstyle{remark}
 \newtheorem{remark}[theorem]{Remark}
 \newtheorem*{theorem*}{Theorem}
 \newtheorem*{example}{Example}
 \numberwithin{equation}{section}

\title{Weighted Boundedness of a Composition of Paraproducts}
\address{
Mathematics and Science Building\\
University of Missouri\\
810 Rollins St\\
Columbia, MO USA\\
65201}
\email{acg7y@umsystem.edu}

\author{Ana Čolović}
\maketitle

\section{Introduction}
\label{sec1}

Paraproducts first appeared in the work of Bony (in \cite{Bon}), who used them to study pseudodifferential operators. Since then, the term paraproduct has been used to refer to a broad set of operators, which can be thought of as "parts of products" of two functions (see \cite{BenMalNai}). Given their role in decomposing products of functions, dyadic paraproducts $\Pi_b$, with a symbol $b\in L^1_{loc},$ have been central in understanding commutators of singular integral operators $T$ and the multiplication operator $M_b$, denoted by $[M_b,T]$ or $[b, T],$ starting from the work of Calder\'{o}n (in \cite{Cal}), who studied the commutator of the Cauchy transform and the multiplication operator. Establishing the weighted bounds for the paraproducts was later the crucial piece of establishing the sharp weighted bounds for commutators of general Calder\'{o}n-Zygmund operators on $L^p(w)$ spaces, where $1<p<\infty,$ and $w$ is a Muckenhoupt weight in the work of Hytonen (in \cite{Hyt}).

The study of compositions of dyadic paraproducts was initiated by Pott, Reguera, Sawyer and Wick (in \cite{PotRegSawWic}), who formulated a dyadic version of Sarason's conjecture concerning Toeplitz operators, and used compositions of dyadic paraproducts as discrete models for Toeplitz operators.

One of the operators studied by Pott, Reguera, Sawyer and Wick was the operator $\Pi_b^*\Pi_d,$ a composition of an adjoint of a dyadic paraproduct with a symbol $b,$ and a paraproduct with a symbol $d.$

In particular, in \cite{PotRegSawWic}, Pott, Reguera, Sawyer and Wick established the conditions for the $L^2(\mathbb{R})$ boundedness of the operator $\Pi_b^*\Pi_d.$  In \cite{HolFra}, Holmes and Fragkiadaki continued the study of the operator $\Pi_b^*\Pi_d,$ and established weighed bounds for the operator $\Pi_b^* \Pi_d$, showing that sparse domination yields an upper bound in $L^p(w)$ spaces, where $1<p<\infty,$ and $w$ belongs to the $A_p$ class of weights. In this paper, we reconcile the two bounds, and show that one can find a sparse bound for the operator $\Pi_b^*\Pi_d$ that matches the bound from the $L^2(\mathbb{R})$ setting of Pott, Reguera, Sawyer and Wick. Furthermore, we establish the lower bounds when the operator acts on the $L^2(w)$ space, and $w\in A_2.$ 

 Here, we use $\mathcal{D}$ to refer to the standard dyadic grid of intervals on $\mathbb{R}$ and $\{h_I\}_{I\in \mathcal{D}}$ to the Haar basis of wavelets, where
\begin{equation} \label{E:Haar}
h_I:=\frac{1}{\sqrt{|I|}}\left(\mathbf{1}_{I^{+}}-\mathbf{1}_{I^-}\right),
\end{equation}
and $I^+,$ and $I^{-}$ denote the left and right dyadic children of $I,$ respectively.

Given a locally integrable function $a,$ we define a sequence $\{a_I\}_{I\in \mathcal{D}},$ where, for each $I\in \mathcal{D},$ $a_I=\langle a,h_I \rangle.$ We identify locally integrable functions $a$ and the sequences of their Haar coefficients.

 To discuss the boundedness of operators $\Pi_b^*\Pi_d,$ we introduce two sequences. Given a sequence of complex numbers, $a=\{a_I\}_{I\in \mathcal{D}},$ the sweep of a sequence is defined as:
\[
\widehat{S}(a)=\left\{\sum_{J\subsetneq I} a_J h_I (J)\right\}_{I\in\mathcal{D}},
\]
where $h_I(J)$ is the constant value of the function $h_I$ on the interval $J,$ when $J\subsetneq I.$
The sequence $E(a)$ is defined as:
\[
E(a)=\left\{\frac{1}{|I|}\sum_{J\subset I} a_J \right\}_{I\in \mathcal{D}}.
\]

When discussing paraproducts, we refer to sequences $b=\{b_I\}_{I\in \mathcal{D}},$ and  $d=\{d_I\}_{I\in \mathcal{D}},$ as symbols of the operator. 
The Schur product of $b$ and $d$ is the sequence $\{b_I\, d_I\}_{I\in \mathcal{D}}$ denoted by $b\circ d. $ Given a sequence $\{a_I\},$ we define two norms:
\begin{align}
    \|a\|_{CM}&:=\sqrt{\sup_{I\in \mathcal{D}} \frac{1}{|I|} \sum_{J\subset I} \left|a_J\right|^2} \label{E:CM-norm}\\
    \|a\|_{l^{\infty}}&:=\sup_{I\in \mathcal{D}}|a_I|.
\end{align}
Pott, Reguera, Sawyer and Wick proved the following theorem.
\begin{theorem}[Pott, Reguera, Sawyer and Wick - 2016]
    Let $b=\{b_I\}_{I\in \mathcal{D}},$ and $d=\{d_I\}_{I\in \mathcal{D}}$ be sequences of complex numbers. The operator $\Pi_b^*\Pi_d$ is bounded on $L^2(\mathbb{R})$ if and only if 
   $\|\widehat{S}(b\circ d)\|_{CM}<\infty$ and $\|E(b\circ d)\|_{l^{\infty}}<\infty.$ Moreover, 
   \[
   \|\Pi_b^*\Pi_d\|_{L^2(\mathbb{R})}\approx \|\widehat{S}(b\circ d)\|_{CM}+\|E(b\circ d)\|_{l^{\infty}}.
   \]
\end{theorem}

Given the use of paraproducts in establishing weighted inequalities for commutators of the Calderon-Zygmund operators (in \cite{Hyt}), in this paper we aim to explore the question of weighted boundedness of the compositions of paraproducts.

As discussed earlier, the weighted theory of these operators has first been studied by Fragkiadaki and Holmes in \cite{HolFra}, who considered the operator $\Pi_b^*,$ and a dyadic paraproduct $\Pi_d$ with symbols $b,d \in L^1_{loc}$ on weighted spaces $L^p(w),$ where $w\in A_p,$ belong to a Muckenhoupt class of weights, with $1<p<\infty.$  They established upper bounds for the operator $\Pi_b^*\Pi_d$, through sparse domination.

 Weights are locally integrable, non-negative functions and for a weight $w$, the space $L^p(w),$ where $1<p<\infty,$ is defined as the space of all functions $f$
such that $\|f\|_{L^p(w)}:=\left(\int|f(x)|^pw(x)dx\right)^{\frac{1}{p}}<\infty.$ 

For a non-negative measurable function $w$, $1<p<\infty,$ and $q$ the Hölder conjugate of $p,$ we say $w\in A_p$ if
\begin{equation} \label{D:Muckenhoupt}
[w]_{A_p}=\sup_{I\subset \mathbb{R}} \left(\frac{1}{|I|} \int_{I} w \:dm\right)\left(\frac{1}{|I|}\int_{I}w^{-\frac{q}{p}}\:dm\right)^{\frac{p}{q}}<\infty.
\end{equation}
Here, the supremum is taken over intervals in $\mathbb{R}.$  

For an interval $I\subset \mathbb{R},$ we let
\[ w(I) = \int_I w \:dm, \qquad w_I = \frac{w(I)}{|I|}.
\]

A powerful tool for establishing weighted bounds has been the technique of sparse domination.

 For $0<\eta<1,$ we say that a  collection $\mathcal{S}\subset \mathcal{D}$ is $\eta$-sparse if for every $Q\in \mathcal{S}$ there is a measurable subset $E_Q$ of $Q$ such that the collection of sets $\{E_Q\}_{Q\in \mathcal{S}}$ is pairwise disjoint and for all $Q\in \mathcal{S},$
$|E_Q|\geq\eta |Q|.$

Given a sparse collection $\mathcal{S}$ we can define a sparse operator $\mathcal{A}_{\mathcal{S}}$ by 
\[
\mathcal{A}_{\mathcal{S}}f:= \sum_{Q\in \mathcal{S}}\langle f\rangle_Q \mathbf{1}_{Q}
\]

We say that an operator $T$ is dominated by a sparse operator $\mathcal{A}_{\mathcal{S}}$ if either 
$|\langle Tf,g\rangle_{L^2(\mathbb{R})}| \lesssim |\langle \mathcal{A}_{\mathcal{S}}f,g \rangle_{L^2(\mathbb{R})}|$
or $|Tf(x)|\lesssim |\mathcal{A}_{\mathcal{S}}f(x)|, x-  a.e.,$ for locally integrable functions $f$ and $g$, where implied constants do not depend on $f$ and $g.$ 

In \cite{HolFra}, Holmes and Fragkiadaki, showed that the operator $\Pi_b^*\Pi_d$ is dominated by a sparse operator.

\begin{theorem}{(Holmes, Fragkiadaki (2022))} \label{T:HolFra} There is a constant C, such that for all weights $w \in A_p,$ $BMO$ functions $b,d, $ and $1<p<\infty,$
\[
\|\Pi_b^*\Pi_d\|_{L^p(w)} \lesssim \|b\|_{BMO}\|d\|_{BMO}  \|\mathcal{A}_{\mathcal{S}}\|_{ L^p(w)}.
\]
Therefore, 
\[\|\Pi_b^*\Pi_d\|_{L^p(w)} \lesssim \|b\|_{BMO}\|d\|_{BMO} [w]_{A_p}^{\max\{1,\frac{1}{p-1}\}}.\]

\end{theorem}
Here, the dyadic $BMO$ norm of a function $b$ is defined as
\begin{equation}\label{E:Dyadic-BMO}
\|b\|_{BMO^{\mathcal{D}}}:=\sup_{I\in \mathcal{D}}\frac{1}{|I|}\int_{I}|b-\langle b\rangle_{I}|dx,
\end{equation}
where, for $I\in \mathcal{D,}$
$\langle b\rangle_I=\frac{1}{|I|}\int_{I}b\, dx.$

We ask whether the bound of Pott, Reguera, Sawyer and Wick can be extended to the upper bound on $L^p(\mathbb{R},w)$, for Muckenhoupt weights $w,$. We establish a sparse domination result for the operator $\Pi_b^*\Pi_d,$ with the bounds from the $L^2(\mathbb{R})$ setting. 
In doing so we combine the two perspectives on the operator $\Pi_b^*\Pi_d$, from Fragkiadaki and Holmes, and Pott, Reguera, Sawyer and Wick and prove the following theorem. 

\begin{theorem} \label{T:Sparse-improved}
    Let $w\in A_p.$ Let $b=\{b_I\}_{I\in \mathcal{D}}, d=\{d_I\}_{I\in \mathcal{D}}$ and suppose that
$
\left\|\widehat{S}(b\circ d)\right\|_{CM}<\infty,
$
and $\left\|E(b\circ d)\right\|_{l^{\infty}} < \infty.$ Then the operator $\Pi_b^*\Pi_d$ is bounded on $L^p(w),$ and there is a constant $C$ such that 
\begin{align*} 
&\| \Pi^{*}_{b} \Pi_d\|_{L^p(w)} \\
&\leq C \left(\left\|\widehat{S}(b\circ d)\right\|_{CM} + \left\|E(b\circ d)\right\|_{\ell^{\infty}} \right) \|\mathcal{A}_{\mathcal{S}}\|_{ L^p(w)}.
\end{align*}
\end{theorem}
We provide a comparison of the bounds from Theorem \ref{T:Sparse-improved}, and 
the bounds of Holmes and Fragkiadaki (in Theorem \ref{T:HolFra}), showing that our bounds in $L^p(\mathbb{R})$ are better or equivalent to those of Fragkiadaki and Holmes. 
\begin{proposition} Given two locally integrable functions $a$ and $b,$ such that $a, b \in  BMO^{\mathcal{D}},$ we have 
\begin{equation} \label{E:Norm-Ineq}
\|\widehat{S}(b\circ d)\|_{CM}+\|E(b\circ d)\|_{l^{\infty}}\lesssim \|b\|_{BMO}\|d\|_{BMO},
\end{equation}
where $b\circ d=\{\langle b,h_I\rangle \langle d, h_I\rangle\}_{I\in \mathcal{D}}.$
\end{proposition}
If we let $b_I=b_{I_0}$ for $I=I_0$ and $0$ otherwise, and $d_I=d_{J_0}$ for $I=J_0$ and $0$ otherwise, for disjoint dyadic intervals $I_0$ and $J_0$, one can calculate that the norm on the left-hand side of the inequality \eqref{E:Norm-Ineq} is strictly smaller than its right-hand side. 

In the case when $p=2,$ we ask whether the lower bound of the operator $\Pi_b^*\Pi_d$ on $L^2(w),$ for $w\in A_2,$ from Corollary 2.2.  of \cite{PotRegSawWic} can be extended to more general, weighted setting, and prove the following theorem. 
\begin{theorem} \label{T:Lower-bound-11} Let $b=\{b_I\}_{I\in\mathcal{D}},$ and $d=\{d_I\}_{I\in\mathcal{D}},$ and let $w\in A_2.$ Then 
\[
\left\|\widehat{S}(b\circ d)\right\|_{CM}+\left\|E(b\circ d)\right\|_{\ell^{\infty}}  \lesssim \|\Pi^*_b \Pi_d\|_{ L^2(w)},
\]
where the implied constant depends on the $A_2$ characteristic of the weight $w.$
\end{theorem}

In Section \ref{S:Prelim} we offer necessary preliminaries to the proof of Theorem \ref{T:Sparse-improved}, and Theorem \ref{T:Lower-bound-11}, and prove relevant lemmas. The proofs are then provided in Section \ref{S:Proofs}.

\section{Preliminaries}\label{S:Prelim}

To begin, we discuss the ambient function spaces, and their properties.
\subsection{\texorpdfstring{$BMO$, $BMO(w)$ and  Weighted Square Function}{}}
\label{SS:weighted-inequalities}

We first discuss the $BMO$ and weighted $BMO$ spaces, and recall two results in the theory of weighted inequalities, that we use in proving Theorem \ref{T:Lower-bound-11}.

We recall the definition of the $BMO^{\mathcal{D}}$ functions.

\begin{definition}[$BMO^{\mathcal{D}}$]
   A function $b$ belongs to $BMO^{\mathcal{D}}$ if 
\[
\|b\|_{BMO^{\mathcal{D}}}:=\sup_{I\in \mathcal{D}}\frac{1}{|I|}\int_{I}|b-\langle b\rangle_I|dx,
\]
where, for $I\in \mathcal{D},$ $\langle b \rangle_I=\frac{1}{|I|}\int_Ib \, dx.$ 
\end{definition}
 The supremum above is taken over all dyadic intervals. It is also possible to define a norm with a supremum over all intervals $I\subset \mathbb{R}.$ In this case the function is said to belong to $BMO.$ The space $BMO^{\mathcal{D}}$ is contained in $BMO,$ but the two are not equivalent (see \cite{GarJon}).

It is a well-known consequence of John-Nirenberg's theorem (see \cite{Gar}) that one can replace the $L^1$ averages in the definition above with $L^2$ averages.
Hence, 
\begin{equation} \label{E:JN2}
\|b\|_{BMO^{\mathcal{D}}}\approx \sup_{I\in \mathcal{D}}\left(\frac{1}{|I|}\int_{I}|b-\langle b\rangle_I|^2\,dx\right)^{\frac{1}{2}}.
\end{equation}  
Haar wavelets defined in the equation \eqref{E:Haar} are an orthonormal basis in $L^2(\mathbb{R}),$ as shown by Haar in his PhD thesis (in \cite{Haa}).

For $I\in \mathcal{D},$ expanding the function $\frac{\mathbf{1}_{I}}{|I|}$ in the Haar wavelet basis, we get that 
\begin{equation} \label{E:Average}
\langle b\rangle_I=\left\langle b, \frac{\mathbf{1}_{I}}{|I|}\right\rangle=\sum_{J: I\subsetneq J} \langle b, h_J\rangle h_J(I),
\end{equation}
where, since for $J\subsetneq I,$ $h_J$ is constant on $I,$ $h_J(I)$ denotes the value of the function $h_J$ on $I.$

Using the equation (\ref{E:Average}), one can show that, for each $I\in \mathcal{D},$
\begin{equation}\label{E:Dyadic-BMO-sum}
\frac{1}{|I|}\int_{I}|b-\langle b\rangle_I|^2 \,dx =\sum_{J:I\subsetneq J}|\langle b,h_J\rangle|^2.
\end{equation}
It then follows that for a sequence $b=\{\langle b,h_I\rangle\}_{I\in \mathcal{D}}=\{b_I\},$ 
\begin{equation*}
   \|b\|_{BMO^{\mathcal{D}}}\approx \|b\|_{CM}. 
\end{equation*}

Given a weight $w\in A_2,$ one can define the space $BMO^{\mathcal{D}}(w)$ as the space of locally integrable functions $b$
such that 
\[
\|b\|_{BMO^{\mathcal{D}}}:=\sup_{I\in \mathcal{D}}\frac{1}{w(I)}\int_{I}|b-\langle b\rangle_I|w\,dx<\infty
\]
In the case when instead of dyadic, we consider taking the supremum over all possible intervals $I\subset \mathbb{R},$ Muckenhoupt and Wheeden proved in \cite{MucWhe} that $\|b\|_{BMO}\approx \|b\|_{BMO(w)}.$

 For the dyadic $BMO^{\mathcal{D}},$ and $BMO^{\mathcal{D}}(w),$ spaces, Holmes, Lacey and Wick, in \cite{HolLacWic} showed that the same equivalence of norms holds. From Theorem 4.1. in \cite{HolLacWic}, it follows that
\begin{equation} \label{E:BMO-Weighted-Equivalence}
    \|b\|_{BMO^{\mathcal{D}}}\approx\|b\|_{BMO^{\mathcal{D}}(w)}.
\end{equation}

Since, in Theorem \ref{T:Lower-bound-11}, we establish lower bounds for the operator $\Pi_b^*\Pi_d,$ in $L^2(w),$ where $w\in A_2,$ we recall a result in the theory of $L^2(w)$ spaces for $w\in A_2,$ due to Petermichl and Pott (in \cite{PetPot}). Their result allows one to estimate the $L^2(w)$ norm of a function by understanding its $L^2(\mathbb{R})$ coefficients. 

We recall the definition of the dyadic square function on $L^2(\mathbb{R}).$
\begin{definition}[Dyadic Square Function] For a locally integrable function $f$, the square function $S$ is given by
\[
Sf(t):=\left(\sum_{I\in \mathcal{D}}|\langle f,h_I\rangle_{L^2(\mathbb{R})}|^2\frac{\mathbf{1}_I(t)}{|I|}\right)^{\frac{1}{2}}.
\]
    
\end{definition}
In \cite{PetPot}, Petermichl and Pott, gave the upper and lower bound of the square function on $L^2(w).$ They proved the following result, which combines Theorems 3.1. and Corollary 3.2. from their work.
\begin{theorem}(Petermichl, Pott - 2001)\label{T:Petermichl-Pott}
    Let $w\in A_2$. For all the functions $f\in L^2(w),$
\begin{equation} \label{E:Square-weighted-upper}
\|f\|_{L^2(w)}\lesssim [w]^{\frac{1}{2}}_{A_2}\|Sf\|_{L^2(w)}=[w]^{\frac{1}{2}}_{A_2}\left(\sum_{I\in \mathcal{D}}|\langle f, h_I\rangle_{L^2(\mathbb{R})}|^2\langle w\rangle_I\right)^{\frac{1}{2}}.
\end{equation}
Furthermore,
\begin{equation} \label{E:Square-weighted-lower}
\|Sf\|_{L^2(w)}=\left(\sum_{I\in \mathcal{D}}|\langle f, h_I\rangle_{L^2(\mathbb{R})}|^2\langle w\rangle_I\right)^{\frac{1}{2}} \lesssim [w]_{A_2}\|f\|_{L^2(w)}.
\end{equation}
\end{theorem}

\subsection{Decomposition of the (1,0,0,1) operator and the sparse bounds of its pieces}
In \cite{PotSmi}, Pott and Smith prove that the operator $\Pi_b \circ \Pi_d^*$ can be decomposed  into three operators: a paraproduct, its adjoint, and a martingale operator.

\begin{proposition}[Pott-Smith Identity] \label{P:Pott-Smith} For symbols $b = \{b_I\}_{I\in \mathcal{D}},$ and $d= \{d_I\}_{I\in\mathcal{D}},$ 
\[
\Pi^*_b \Pi_d  = \Pi_{\widehat{S}(b\circ d)} + \Pi^*_{\widehat{S}(b \circ d)} + T_{E(b\circ d)},
\]
where $T_{E(b\circ d)}$ is the martingale operator given by 
\[
T_{E(b\circ d)} f = \sum_{I\in \mathcal{D}} (E(b\circ d))_{I} (f,h_I) h_I,
\]
where $f$ is a locally integrable function.
    
\end{proposition}

In order to understand the sparse bound for $\Pi_b^*\Pi_d,$ we first discuss the sparse bounds for each of its summands in the Pott-Smith identity. These are results well-known in the literature. We provide references for interested readers.

The argument for the bilinear sparse domination of the paraproduct operator can be adapted from the argument present in Appendix B on page 25 of \cite{HolFra}. Since the proof follows the method of \cite{HolFra} we leave the details to the interested reader. 

\begin{lemma}[Sparse bound of a paraproduct] \label{L:Sparse-paraproduct}
Let $b=\{b_I\}_{I\in \mathcal{D}},$ and for each $I\in \mathcal{D},$ $d_I\in \mathbb{C}. $ For each locally integrable functions $f_1$ and $f_2$ from $\mathbb{R}$ to $\mathbb{C},$ one can find a sparse collection of sets $\mathcal{S},$ such that 
\begin{equation}\label{E:Sparse-paraproduct}
|\langle \Pi_b f_1, f_2\rangle_{L^2(\mathbb{R})} | \lesssim \|b\|_{BMO^{\mathcal{D}}}\sum_{I\in \mathcal{S}}|I|\langle |f_1|\rangle_I \langle|f_2|\rangle_I.
\end{equation}
\end{lemma}

\begin{remark}\label{R:adjoints} Observe that, by the definition of an adjoint, Lemma \ref{L:Sparse-paraproduct} holds if we replace $\Pi_b$ with the operator $\Pi^*_b.$
\end{remark}

To find a bilinear sparse form that dominates the martingale operator, we recall the pointwise domination result by Lacey (in \cite{Lac}, Theorem 2.4.) from which the bilinear sparse domination result follows. For a dyadic martingale operator $T_{\epsilon},$ Lacey's result assumes that $\|\epsilon\|_{l^{\infty}}\leq 1.$ For a general dyadic martingale operator $T_{\epsilon},$ with bounded norm, the symbol can be normalized to obtain the following result. 

\begin{lemma}[Sparse bound of a dyadic martingale - Lacey (2017) ] \label{L:Sparse-martingale-Lacey}
Let $\epsilon=\{\epsilon_I\}_{I\in \mathcal{D}},$ and for each $I\in \mathcal{D},$ $\epsilon_I\in \mathbb{C}. $ For all the functions $f\in L^1,$ supported on a dyadic interval $I_0,$ there is a sparse collection of intervals $\mathcal{S}$ such that 
\begin{equation} \label{E:Pointwise-sparse-martingale}
|T_{\epsilon}f|\mathbf{1}_{I_0}(x)\lesssim  \|\epsilon\|_{l^{\infty}} \sum_{I\in \mathcal{S}}\langle|f|\rangle_I \mathbf{1}_I.
\end{equation}
\end{lemma}

As a corollary of the above point-wise result, one can derive the following bi-linear sparse domination result.

\begin{corollary} \label{C:sparse-blinear-martingale} Let $f_1,f_2$ be two locally integrable functions. Then there exists a sparse collection of dyadic intervals $\mathcal{S}$ such that  
\begin{equation} \label{E:Sparse-bilinear-martingale}
    |\langle T_{\epsilon}f_1,f_2\rangle_{L^2(\mathbb{R})}|\lesssim \|\epsilon\|_{l^{\infty
    }} \sum_{I \in \mathcal{S}} |I |\langle |f_1|\rangle_I \langle |f_2|\rangle_I.
\end{equation}
\end{corollary}

The proof of the Corollary \ref{C:sparse-blinear-martingale} follows by finding a dyadic interval that contains the supports of functions $f_1$ and $f_2$ and applying the Lemma \ref{L:Sparse-martingale-Lacey}.

\subsection{Multivariable setting}
In order to prove the sparse bound for $\Pi_b^*\Pi_d,$ we recall a multilinear sparse domination result proved by Culiuc, Di Plinio, and Ou in \cite{CulDiPOu}. The Remark \ref{R:N-sparse} discusses the way we integrate Pott-Smith identity and the multilinear result of \cite{CulDiPOu} to prove a sparse bound for the operator $\Pi_b^*\Pi_d.$

We first introduce the multilinear setting of \cite{CulDiPOu}, also presented in \cite{CulDiPOu}.

Let $n\geq1$, $\vec{p}\in (0,\infty)^{n+1},$ $g^1,...,g^n$ be functions in $L^{\infty}_{0}(\mathbb{R}^n, \mathbb{C})$. We let $T$ be an operator that takes $L_{0}^{\infty}(\mathbb{R}^d, \mathbb{C})^n$ into the space $L^{\infty}_{0}(\mathbb{R}^d,\mathbb{C}),$ and we assume that $T$ is $n-$sublinear, i.e. that it is sublinear in each of its $n-$components.  

We say that the operator $T$  has a sparse bound if there exists a constant $C>0$ and a sparse collection $\mathcal{S}$ of dyadic rectangles such that
\[
|\langle T(g^1,...,g^n), g^{n+1}\rangle | \leq C \sum_{Q\in \mathcal{S}} |Q| \prod_{j=1}^{n+1} \langle g^j\rangle_{p_j,Q}.
\]
The least such constant $C>0$ is denoted by $\|T\|_{\vec{p}}.$ The sparse collection of sets above depends on the $n+1$-tuple, $\vec{g}=(g^1,...,g^{n+1}).$ Here, $T$ is an $n-$sublinear operator that takes $n-$tuples of locally integrable functions into the space of locally integrable functions, i.e. it takes vectors to scalars. We now define a vector-valued operator that combines several $n-$sublinear operators.

Let $\mathbf{T}=\{T_1,...,T_N\}$ be an $n$-tuple of $n$ sublinear operators. If $\mathbf{f}=(\mathbf{f}^1, ..., \mathbf{f}^n)$ is an n-tuple of locally integrable functions from $\mathbb{R}^n$ to $\mathbb{C}^N,$ and $\mathbf{f}^{n+1}$ is a locally integrable function with the same domain and range, we have 
\[
\langle \mathbf{T}(\mathbf{f}^1, ..., \mathbf{f}^n),\mathbf{f}^{n+1} \rangle:= \sum_{k=1}^{N}\langle T_k(f^1_k,...,f^n_k),f^{n+1}_k\rangle.
\]
The vector-valued operator above acts as a sum of its components.

Let $(r_1,...,r_{n+1})$ be an $n+1$-tuple such that
\begin{equation} \label{E:Banach-Holder}
1\leq r_j\leq \infty, \quad j=1,...,n+1, \qquad r:=\frac{r_{n+1}}{r_{n+1}-1}:= \frac{1}{\sum_{j=1}^{n}\frac{1}{r_j}}.
\end{equation}
We say that $\mathbf{T}$ has a $(\vec{p},r)-$sparse bound if there exist a constant $C>0$ and for each $(\mathbf{f}^1, ..., \mathbf{f}^n) \in L^{\infty}_{0}(\mathbb{R^d},\mathbb{C}^N)$ a sparse collection $\mathcal{S}$ such that 
\[
|\mathbf{T}(\mathbf{f}^1, ..., \mathbf{f}^n),\mathbf{f}^{n+1} \rangle_{L^2(\mathbb{R})}| \leq C \sum_{Q\in \mathcal{S}}|Q|\prod_{j=1}^{n+1}\langle \|\mathbf{f}\|_{l^{r_j}}\rangle_{p_j,Q}.
\]
If $\mathbf{T}$ has a sparse bound, we denote the least such constant $C$ by $\|\mathbf{T}\|_{(\vec{p},\vec{r})}.$

Note that the above norms $\|T\|_{(p,r)}$ and $\|\mathbf{T}\|_{(\vec{p},\vec{r})}$ refer to the sparse bounds of operators, and not apriori to the operator norms of $T$ and $\mathbf{T}$. With this in mind, we present the theorem of Culiuc, Di Plinio and Ou (see Theorem A, in \cite{CulDiPOu}).

\begin{theorem}(Culiuc, Di Plinio, Ou - 2017)
\label{T:linear-CDO} Suppose that $\vec{p}\in[1,\infty)^{n+1},$ and let $\vec{r}$ be as in (\ref{E:Banach-Holder}). Suppose $r_j<p_j$. If $T_1,...,T_N$ are $n-$sublinear operators defined as above, then 

\[
\|\{T_1,...,T_N\}\|_{(\vec{p},\vec{r})} \lesssim \sup_{k=1,...,N} \|T_k\|_{\vec{p}},
\]
where the constant above depends on $\vec{p}, \vec{r}$ and the dimension $d.$
\end{theorem}

The theorem above states that if $n-$sublinear operators have sparse bounds, then so does the vector form that combines them, and moreover, the sparse bound of the vector form is no worse than the biggest among the sparse bounds of the summands. 

\begin{remark} \label{R:N-sparse} The vector-valued multilinear result of Culiuc, Di Plinio and Ou has an application to the sum of linear operators. Suppose $T_1,...T_N$ are linear operators such that for each $1\leq i \leq N,$ $T_i$ has a $(1,1)$ sparse bound. In other words, for each $1\leq i\leq N,$ there exists $C_i>0$ such that for  $(f^{1}_i,f^{2}_i)\in L^{\infty}_{0}(\mathbb{R}^d,\mathbb{C})^2$ one can find a sparse collection of sets $\mathcal{S}_i$ such that 
\[
|\langle T_if_1, f_2\rangle_{L^2(\mathbb{R}^d)}| \leq C_i \sum_{Q\in \mathcal{S}_i} |Q|\langle |f_1|\rangle_Q \langle |f_2|\rangle_Q.
\]
Let $\vec{p}=(1,1),$ and $\vec{r}=(N,N),$ and assume that $N> 1.$
By Theorem \ref{T:linear-CDO}, it follows that $\|\{T_1, ..., T_N\}\|_{((1,1),(N,N)} \lesssim \sup_{j=1,...,N}\|T_j\|_{(1,1)}.$ 
In other words, there exists $C$ such that for each $\vec{\mathbf{f}}=\left\langle (f^1_1,...,f^1_N), (f_1^2,...,f_N^2)\right\rangle\in L^{\infty}_0(\mathbb{R}^d,\mathbb{C}^N)^2,$ we can find a sparse collection $\mathcal{S}$ depending on $\vec{\mathbf{f}},$ such that
\begin{align*}
&|\langle \{T_1,...,T_N\}(f_1^1,....,f_N^1),(f_1^2,...,f_N^2)\rangle|\\&= \left|\sum_{i=1}^{N}\langle T_if_i^1,f_i^2\rangle\right|\\ &\leq C\sum_{Q\in \mathcal{S}}|Q|\langle \|(f_1^1,...,f_N^1)\|_{l^N}\rangle_{Q}\langle \|(f_1^2,...,f_N^2)\|_{l^N}\rangle_{Q}.
\end{align*}
In particular, if $f,g\in L^{\infty}_0(\mathbb{R}^d,\mathbb{C}),$ and we let $f_i^1=f,$ and $f_i^2=g$ for $1\leq i\leq N,$ in the statement above, Theorem \ref{T:linear-CDO} states that there exists a sparse collection $\mathcal{S}$ such that
\[
\left|\sum_{i=1}^{N}\langle T_if,g\rangle\right| = \left|\left\langle\sum_{i=1}^{N}T_if,g\right \rangle\right| \leq C N^{\frac{2}{N}}\sum_{Q\in \mathcal{S}}|Q|\langle |f|\rangle_{Q}\langle |g|\rangle_{Q}.
\]
\end{remark}

\section{Proof of Theorem \ref{T:Sparse-improved}\label{S:Proofs}, Comparison of Bounds and the Proof of Theorem \ref{T:Lower-bound-11}}

Using Theorem \ref{T:linear-CDO} and its particular instance discussed in Remark \ref{R:N-sparse}, we give a sufficient condition for the operator $\Pi^{(1,0)}_b \circ \Pi^{(0,1)}_d = \Pi^{(1,1)}_{b\circ d}$ to be bounded between weighted spaces, and in the following subsection we provide the proof of its lower bound. 
\subsection{Proof of Theorem \ref{T:Sparse-improved}}
\begin{proof}[Proof of Theorem \ref{T:Sparse-improved}]

By Pott-Smith identity (\ref{P:Pott-Smith}), for symbols $b=\{b_I\}_{I \in \mathcal{D}},$ $d=\{d_I\}_{I\in \mathcal{D}},$ the operator $\Pi_b^{*}\Pi_d$ decomposes into a sum of three operators: 
\[
\Pi_b^{*}\Pi_d=  \Pi_{\widehat{S}(b\circ d)} + \Pi^*_{\widehat{S}(b \circ d)} + T_{E(b\circ d)}.
\]
By Lemma \ref{L:Sparse-paraproduct}, Remark \ref{R:adjoints} and Corollary \ref{C:sparse-blinear-martingale} it follows that for each of the operators $T_1= \Pi_{\widehat{S}(b\circ d)},$ $T_2=\Pi^*_{\widehat{S}(b\circ d)}$ and $T_3= T_{E(b\circ d)},$ given locally integrable functions $f,g$ one can find sparse collections $\mathcal{S}_i,$ for $i=1,2,3$ 
such that 
\[
|\langle T_i f,g \rangle|\lesssim C_i \sum_{I\in \mathcal{S}_i} \langle |f|\rangle_I\langle |g|\rangle_I |I|,
\]
where $C_1=C_2=\|\widehat{S}(a\circ b)\|_{BMO^{\mathcal{D}}}, $ and $C_3 = \|E(a\circ b)\|_{l^{\infty}}.$ Therefore, taking $N=3,$ in the Remark \ref{R:N-sparse}, it follows that for each $f,g \in L^{\infty}_0(\mathbb{R},\mathbb{C})$ there exists a sparse collection $\mathcal{S}$ such that 

\[
|\langle \Pi^*_b \Pi_d f, g\rangle_{L^2(\mathbb{R})}| \lesssim \max \left\{\left\|\widehat{S}(b\circ d)\right\|_{BMO^{\mathcal{D}}},\left\|E(b\circ d)\right\|_{\ell^{\infty}}\right\}  \sum_{Q\in \mathcal{S}} |Q| \langle |f|\rangle_Q  \langle |g| \rangle_Q.
\]
 Locally finite functions are dense in $L^2(\mathbb{R})$, so the proof of Theorem \ref{T:Sparse-improved} follows.\end{proof}
\subsection{Comparison of Bounds}
In \cite{PotRegSawWic}, Pott, Reguera, Sawyer and Wick proved that for symbols $b=\{b_I\}_{I\in \mathcal{D}},$ and $d=\{d_I\}_{I\in \mathcal{D}},$
\[
\|\Pi^*_b \Pi_d\|_{ L^2(\mathbb{R})}\ \simeq    \left\|\widehat{S}(b\circ d)\right\|_{BMO^{\mathcal{D}}}+\left\|E(b\circ d)\right\|_{\ell^{\infty}}.
\]
We also know that \begin{align*}
\|\Pi^*_b \Pi_d\|_{L^2(\mathbb{R})} &\leq \|\Pi^*_b\|_{L^2(\mathbb{R})} \|\Pi_d\|_{L^2(\mathbb{R})} \\
&\lesssim \|b\|_{BMO^{\mathcal{D}}} \|d\|_{BMO^{\mathcal{D}}}.
\end{align*}
Therefore, 
$\left\|\widehat{S}(b\circ d)\right\|_{BMO^{\mathcal{D}}}+\left\|E(b\circ d)\right\|_{\ell^{\infty}} \lesssim \|b\|_{BMO} \|d\|_{BMO}. $
We now provide a direct proof of this fact. 

\begin{proposition} For symbols $b=\{b_I\}_{I\in \mathcal{D}},$ and $d=\{d_I\}_{I\in \mathcal{D}},$
the following holds:
\[
\left\|\widehat{S}(b\circ d)\right\|_{BMO^{\mathcal{D}}}+\left\|E(b\circ d)\right\|_{\ell^{\infty}} \lesssim \|b\|_{BMO} \|d\|_{BMO}. 
\]
    
\end{proposition}

\begin{proof}
    Let $I\in \mathcal{D}.$ By Cauchy-Schwarz inequality, it follows that:
    \begin{align*}
    \left|E(b\circ d)_I\right|=
    \frac{1}{|I|}\left|\sum_{J\subset I} b_J d_J\right| 
    &\leq \left(\frac{1}{|I|}\sum_{J\subset I}|b_J|^2\right)^{\frac{1}{2}} \left(\frac{1}{|I|}\sum_{J\subset I}|d_J|^2\right)^{\frac{1}{2}} \\
    &\leq \|b\|_{BMO} \|d\|_{BMO}.
  \end{align*}
  Therefore, $\left\|E(b\circ d)\right\|_{\ell^{\infty}} \leq \|b\|_{BMO} \|d\|_{BMO}. $
  For $I\in \mathcal{D},$ we consider
  \[
  \left(\frac{1}{|I|} \sum_{J\subset I} \left|(\widehat{S}(b\circ d))_J\right|^2\right)^{\frac{1}{2}}
  =\left(\frac{1}{|I|} \sum_{J\subset I} \left|\sum_{K\subsetneq J} b_K d_K h_J(K)\right|^2\right)^{\frac{1}{2}}.
  \] 
  The quantity $\left(\sum_{J\subset I} \left|(\widehat{S}(b\circ d))_J\right|^2\right)^{\frac{1}{2}}$ is the $l^2$ norm of a sequence $\{(\widehat{S}(b\circ d))_{J}\}_{J\subset I}.$ To estimate this norm we use $l^2$ duality. Suppose $c=\{c_J\}_{J\subset I}$ is a sequence in $l^2(\mathcal{D}).$ In order to simplify the notation in the proof below, we will identify the sequence $c=\{c_I\}_{I\in \mathcal{D}}$ with the function given by $\sum_{I\in \mathcal{D}}c_I h_I.$
  Now, 
 \begin{align*}
  \left|\langle \widehat{S}^{I}(b\circ d), c \rangle_{l^2(\mathcal{D})}\right| &= \left|\sum_{J\subset I} (\widehat{S}(b\circ d))_J c_J\right| =
  \left|\sum_{J\subset I} \sum_{K\subset J} b_K d_K h_J(K)  c_J\right|.
   \end{align*}
  By Fubini's theorem, 
  \begin{align*}
  \left|\langle\widehat{S}^{I}(b\circ d), c \rangle_{l^2(\mathcal{D})}\right| &= \left|\sum_{K\subset I} b_K d_K \sum_{K\subsetneq I} c_Jh_J(K)\right| \\ 
  &= \left| \sum_{K\subset J} b_K d_K \langle c \rangle_{K} \right|= \left|\langle b, \sum_{K\subset I}d_K\langle c \rangle_K h_K\rangle_{L^2(\mathbb{R})}\right| \\
  &\lesssim \|b\|_{BMO} \|S_{\mathcal{D}}\phi\|_{L^1},
   \end{align*}
where the last inequality follows by $H^1-BMO$ duality, and where the function $\phi:=\sum_{K\subset I}d_K\langle c \rangle_K h_K.$
Since $S_{\mathcal{D}}\phi \leq M_{\mathcal{D}}c \: S_{\mathcal{D}}d^{I},$  where $d^I=\{d_J\}_{J\subset I}$, and by the definition of $\|d\|_{BMO}$ we have that 
\begin{align*}
    \|S_{\mathcal{D}}\phi\|_{L^1} \leq \|M_{\mathcal{D}}c\|_{L^2}\|Sd^{I}\|_{L^2} \lesssim \|c\|_{L^2} \|d^I\|_{L^2} \leq \|c\|_{L^2} |I|^{\frac{1}{2}} \|d\|_{BMO}.
\end{align*} 
Therefore,  
\[
 \left(\frac{1}{|I|} \sum_{J\subset I} \left|(\widehat{S}(b\circ d))_J\right|^2\right)^{\frac{1}{2}} \lesssim \|b\|_{BMO}\|d\|_{BMO},
\]
which concludes the proof.\end{proof}
\subsection{Proof of Theorem \ref{T:Lower-bound-11}}
In this subsection we prove the Theorem \ref{T:Lower-bound-11} which concerns the lower bound of the operator $\Pi_b^*\Pi_d$ on $L^2(w).$
We begin by recalling the statement of Theorem \ref{T:Lower-bound-11}.
\begin{theorem*}
    Let $b=\{b_I\}_{I\in\mathcal{D}},$ and $d=\{d_I\}_{I\in\mathcal{D}},$ and let $w\in A_2.$ Then 
\[
\left\|\widehat{S}(b\circ d)\right\|_{BMO^{\mathcal{D}}}+\left\|E(b\circ d)\right\|_{\ell^{\infty}}  \lesssim \max\{[w]_{A_2}, [w]_{A_{\infty}}[w]^2_{A_2}\} \|\Pi^*_b \Pi_d\|_{ L^2(w)}.
\]
\end{theorem*}
\begin{proof} 
    We observe that
    \begin{align*}
     &\|E(b\circ d) \|_{\ell^{\infty}}=\sup_{I\in\mathcal{D}}|E(b\circ d)(I)|= \sup_{I\in\mathcal{D}}|\langle \Pi^*_b \Pi_d h_I, h_I\rangle_{L^2(\mathbb{R})}| \\
     &\leq \|\Pi^*_b \Pi_d\|_{ L^2(w)} \sup_{I\in\mathcal{D}}\|h_{I}\|_{L^2(w)} \|h_{I}\|_{L^2(w^{-1})}. 
    \end{align*}
    Since $w\in A_2,$ 
    \[\sup_{I\in\mathcal{D}}\|h_{I}\|_{L^2(w)} \|h_{I}\|_{L^2(w^{-1})}=\sup_{I\in\mathcal{D}}\left(\frac{w(I)}{|I|}\frac{w^{-1}(I)}{|I|}\right)^{\frac{1}{2}}< \infty.\]

Therefore, 
\[
\|E(b\circ d) \|_{\ell^{\infty}} \leq [w]_{A_2}\|\Pi^*_b \Pi_d\|_{ L^2(w)}
\]
Let $I\in \mathcal{D},$ $k\in \mathbb{N}.$ For $K\subset I,$ $|K|> 2^{-k}|I|,$ we define the function
\[
F_{K}^{k}:= \sum_{\substack{J\subset K \\|J|> 2^{-k}|I|}} \overline{\widehat{S}(b\circ d)(J)}h_J.
\]
Since the function $F_{K}^k$ is a finite sum of compactly supported, bounded functions, it follows that $F_{K}^k \in L^2(w).$ 
By a calculation, it follows that
\begin{align} \label{E:K-testing}
&\frac{1}{|\widehat{K}|^{\frac{1}{2}}}\sum_{\substack{J\subset K \\|J|> 2^{-k}|I|}} |\widehat{S}(b\circ d)(J)|^2 = |\langle \Pi^*_b \Pi_d F_k^{K}, h_{\widehat{K}} \rangle_{L^2(\mathbb{R})}| \notag
\\ &\leq \|\Pi^*_b \Pi_d\|_{L^2(w)} \|F_k^K\|_{L^2(w)}\|h_{\widehat{K}}\|_{L^2(w^{-1})}.
\end{align}
By the Theorem \ref{T:Petermichl-Pott} and the definition of $F^{K}_k,$
\begin{align*}
\|F_k^K\|_{L^2(w)} &\lesssim [w]_{A_2}\left(\sum_{J\in\mathcal{D}}|(F^K_k,h_J)|^2 \langle w \rangle_J\right)^{\frac{1}{2}}\\
&= \left(\sum_{\substack{J\subset K \\|J|> 2^{-k}|I|}}|\widehat{S}(b\circ d)(J)|^2 \langle w \rangle_J\right)^{\frac{1}{2}}.
\end{align*}

Since $w\in A_2,$ $\|h_{\widehat{K}}\|_{L^2(w^{-1})} = \left(\langle w^{-1} \rangle _{\widehat{K}}\right)^{\frac{1}{2}} \leq [w]_{A_2} \left( \frac{|\widehat{K}|}{w(\widehat{K})}\right)^{\frac{1}{2}}. $

Therefore, the equation (\ref{E:K-testing}), becomes:
\begin{align*}
&\frac{1}{|K|}\sum_{\substack{J\subset K \\|J|> 2^{-k}|I|}} |\widehat{S}(b\circ d)(J)|^2  
\\ &\lesssim [w]^2_{A_2}\|\Pi^*_b \Pi_d\|_{L^2(w)} \left(\frac{1}{w(K)}\sum_{\substack{J\subset K \\|J|> 2^{-k}|I|}}|\widehat{S}(b\circ d)(J)|^2 \langle w \rangle_J\right)^{\frac{1}{2}}.
\end{align*}

Since $F^{I}_k$ is a finite sum of multiples of Haar functions, it has a finite $BMO(w)$ norm. Using a similar reasoning as in the discussion preceding the equation (\ref{E:Dyadic-BMO-sum}), and Theorem \ref{T:Petermichl-Pott}, it follows that
\[
\|F^{I}_k\|_{BMO(w)} = \max_{\substack{K\subset I \\|K|> 2^{-k}|I|}} \left\{\left(\frac{1}{w(K)}\sum_{\substack{J\subset K \\|J|> 2^{-k}|I|}}|\widehat{S}(b\circ d)(J)|^2 \langle w \rangle_J\right)^{\frac{1}{2}}\right\}.
\]
By the equation \eqref{E:BMO-Weighted-Equivalence}, the dyadic $BMO$ and $BMO(w)$ norms are equivalent. Therefore,  
\begin{align*}
&\frac{1}{|K|}\sum_{\substack{J\subset K \\|J|> 2^{-k}|I|}} |\widehat{S}(b\circ d)(J)|^2  
\\ &\lesssim [w]_{A_{\infty}}[w]^2_{A_2} \|\Pi^*_b \Pi_d\|_{L^2(w)} \|F^I_k\|_{BMO},
\end{align*}
for all $K\subset I,$ $|K|>2^{-k}|I|.$
The quantity $[w]_{A_\infty}$ above arises out of the quantitative relationship between the $BMO$ and the $BMO(w)$ norms established in \cite{Tsu}.
Thus, 
\[
\|F^I_k\|^2_{BMO}\lesssim [w]_{A_\infty}[w]^2_{A_2} \|\Pi^*_b \Pi_d\|_{L^2(w)} \|F^I_k\|_{BMO}.
\]
Hence,  
\[ 
 \left(\frac{1}{|I|}\sum_{\substack{J\subset I \\|J|> 2^{-k}|I|}} |\widehat{S}(b\circ d)(J)|^2 \right)^{\frac{1}{2}}\leq  \| F^I_k\|_{BMO}\lesssim [w]_{A_{\infty}} [w]^2_{A_2} \|\Pi^*_b \Pi_d\|_{L^2(w)}.
\]
Taking a limit as $k\to \infty$, we get that 
\[
\left\|\widehat{S}(b\circ d)\right\|_{BMO^{\mathcal{D}}} \lesssim [w]_{A_\infty}[w]_{A_2}^2\|\Pi^*_b \Pi_d\|_{L^2(w)},
\]
concluding the proof.\end{proof}

\textbf{Acknowledgements.} The author would like thank Brett Wick for many useful conversations regarding this paper.


\begin{thebibliography}{00}


\bibitem{BenMalNai}
Á.~Bényi, D.~Maldonado, and V.~Naibo, What is a paraproduct?,
\textit{Notices Amer. Math. Soc.} \textbf{57} (2010), no.~7, 858–860.

\bibitem{Bon}
J.-M. Bony, Calcul symbolique et propagation des singularités pour les équations aux dérivées partielles non linéaires,
\textit{Ann. Sci. École Norm. Sup.} \textbf{14} (1981), 209–246.

\bibitem{Cal}
A.~P.~Calderón, Commutators of singular integral operators,
\textit{Proc. Nat. Acad. Sci. USA} \textbf{53} (1965), 1092–1099.

\bibitem{CulDiPOu}
A.~Culiuc, F.~Di Plinio, and Y.~Ou, Domination of multilinear singular integrals by positive sparse forms,
\textit{J. London Math. Soc.} \textbf{98} (2018), no.~2, 369–392.

\bibitem{Gar}
J.~B.~Garnett, \textit{Bounded Analytic Functions}, Academic Press, New York, 1981.

\bibitem{GarJon}
J.~B.~Garnett and P.~W.~Jones, BMO from dyadic BMO,
\textit{Pacific J. Math.} \textbf{99} (1982), 351–371.

\bibitem{Haa}
A.~Haar, Zur Theorie der orthogonalen Funktionensysteme,
\textit{Math. Ann.} \textbf{69} (1910), 331–371.

\bibitem{HolFra}
I.~Holmes and V.~Fragkiadaki, Paraproducts, Bloom BMO and Sparse BMO Functions,
\textit{Rev. Mat. Iberoam.} \textbf{39} (2023), no.~6.

\bibitem{HolLacWic}
I.~Holmes, M.~Lacey, and B.~D.~Wick, Commutators in the two weight setting,
\textit{Math. Ann.} \textbf{367} (2017), 51–80.

\bibitem{Hyt}
T.~P.~Hyt\"{o}nen, The sharp weighted bound for general Calderón–Zygmund operators,
\textit{Ann. of Math.} \textbf{175} (2012), 1473–1506.

\bibitem{Lac}
M.~T.~Lacey, An elementary proof of the $A_2$ bound,
\textit{Israel J. Math.} \textbf{217} (2017), 181–195.

\bibitem{MucWhe}
B.~Muckenhoupt and R.~Wheeden, Weighted bounded mean oscillation and the Hilbert transform,
\textit{Studia Math.} \textbf{54} (1975/76), no.~3, 221–237.

\bibitem{PetPot}
S.~Petermichl and S.~Pott, An estimate for weighted Hilbert transform via square functions,
\textit{Trans. Amer. Math. Soc.} \textbf{354} (2002), 1699–1703.

\bibitem{PotRegSawWic}
S.~Pott, M.~C.~Reguera, E.~T.~Sawyer, and B.~D.~Wick, Composition of dyadic paraproducts,
\textit{Adv. Math.} \textbf{298} (2016), 581–611.

\bibitem{PotSmi}
S.~Pott and M.~Smith, Paraproducts and Hankel operators of Schatten class via $p$-John–Nirenberg theorem,
\textit{J. Funct. Anal.} \textbf{217} (2004), 38–78.

\bibitem{Tsu}
Y.~Tsutsui, $A_8$ constants between BMO and weighted BMO,
\textit{Proc. Japan Acad. Ser. A Math. Sci.} \textbf{90} (2014), no.~1, 11–14.





\end{thebibliography}
\end{document}